\documentclass{article}

\usepackage{amsfonts}
\usepackage{amssymb}
\usepackage{amsbsy}
\usepackage{amsmath} 
\usepackage{verbatim}

\newcommand{\hecke}[1]{\mathbf{T}_{#1}}

\newcommand{\zed}{\mathbf{Z}}
\newcommand{\qq}{\mathbf{Q}}
\newcommand{\magma}{{\sc Magma}}
\newcommand{\HECKE}{{\sc hecke}}
\newcommand{\qtwo}{\mathbf{Q_{\rm 2}}}
\newcommand{\ztwo}{\mathbf{Z_{\rm 2}}}

\newcommand{\m}{\mathfrak{m}}
\newcommand{\p}{\mathfrak{i}}

\newcommand{\f}[1]{\mathbf{F_{\rm #1}}}
\newtheorem{thm}{Theorem}
\newtheorem{question}[thm]{Question}
\newcommand{\mystrut}{\rule{0cm}{0.5cm}}

\DeclareMathOperator{\End}{End} 
\DeclareMathOperator{\Frob}{Frob}
\DeclareMathOperator{\Gal}{Gal}
\DeclareMathOperator{\GL}{GL}
\newcommand{\C}{\mathbf{C}}
\newcommand{\F}{\mathbf{F}}
\newcommand{\GQ}{\Gal(\Qbar/\Q)}
\newcommand{\Q}{\mathbf{Q}}
\newcommand{\Qbar}{\overline{\Q}}
\newcommand{\Qtwobar}{\overline{\Q}_2}
\newcommand{\T}{\mathbf{T}}
\newcommand{\Z}{\mathbf{Z}}

\usepackage{epsfig}

\begin{document}

\title{Some non-Gorenstein Hecke algebras attached to spaces of modular
forms.\footnote{This work was partially carried out at the Institut Henri Poincare in Paris; the author is grateful for the hospitality that he received there.}}
\author{L. J. P. Kilford}
\maketitle

%\tableofcontents

\section{Introduction}
 
In this paper we exhibit some examples of non-Gorenstein Hecke algebras, and
hence some modular forms for which mod~2 multiplicity one does not hold.
 
Define $S_2(\Gamma_0(N))$ to be the space of classical cuspidal
modular forms of weight 2, level $N$, and trivial character.
The Hecke algebra~$\hecke{N}$ is defined to be the subring
of $\End(S_2(\Gamma_0(N))$ generated by the Hecke operators
$\left\{T_p:\: p \not |\: N\right\}$ and $\left\{U_q:\:q | N\right\}$.
Let $\m$ be a maximal ideal of~$\hecke{N}$, and let $\ell$ denote the
characteristic of the finite field $\hecke{N}/\m$.
By work of Shimura, one can associate to $\m$ a semi-simple Galois
representation $\rho_{\m}:{\rm Gal}(\overline{\qq}/\qq) \rightarrow
\GL_2(\hecke{N}/\m)$ satisfying
${\rm tr}(\rho({\Frob}_p)) \equiv T_p \: {\rm mod}\: \m$
for all primes $p \nmid N\ell$. We say that~$\m$ is non-Eisenstein if
$\rho_{\m}$ is absolutely irreducible.

As an example, if~$E$ is a (modular) elliptic curve over~$\qq$ of
conductor~$N$,
let $f:=\sum_{n \ge 1} a_n q^n$ be the modular form in~$S_2(\Gamma_0(N))$
associated to~$E$. Associated to $f$ is a minimal prime ideal of $\hecke{N}$;
we say that~$\m$ is associated to~$f$, or to~$E$, if~$\m$ contains this
minimal prime ideal. In this case, the representation associated to~$\m$
is isomorphic to the semisimplification of~$E[\ell]$, where~$\ell$ is
the characteristic of~$\hecke{N}/\m$.

The localisation $(\hecke{N})_{\m}$ of $\hecke{N}$ at a maximal ideal
$\m$ is frequently a Gorenstein ring, and such a phenomenon is related to the
study
of the $\m$-torsion in the Jacobian $J_0(N)$ of $X_0(N)$. For example,
if $\m$ is non-Eisenstein then by the main result of~\cite{BLR},
the $\m$-torsion
in $J_0(N)$ is isomorphic to a direct sum of $d\geq1$ copies of $\rho_{\m}$.
If $d=1$ then one says that the ideal $\m$ satisfies
``mod $\ell$ multiplicity one'', or just ``multiplicity one''.
In this case, the localisation $(\hecke{N})_{\m}$ is known to be Gorenstein.

Multiplicity one is a common phenomenon for maximal
ideals $\m$ of $\hecke{N}$. Let us restrict for the rest of
this paper to the case of non-Eisenstein maximal ideals~$\m$.
The first serious study of these mod $\ell$ multiplicity one questions
is that of Mazur~\cite{mazur}, who proves that if $N=q$ is
prime and the characteristic of the finite field $\hecke{q}/\m$
is not~2 then $\m$ satisfies multiplicity one and hence
${(\hecke{q})}_\m$ is Gorenstein. He also showed that,
if $\hecke{q}/\m$ has characteristic 2 and $T_2\in\m$ then
$(\hecke{q})_{\m}$ is Gorenstein. This work has been generalised
by several authors, to higher weight cases and non-prime level.
Rather than explaining these generalisations in complete generality,
we summarise what their implications are in the case left open by Mazur. 

We fix notation first. From now on, the level $N=q$ is prime,
$\m$ is a non-Eisenstein maximal ideal of $\hecke{q}$,
and the characteristic of $\hecke{q}/\m$
is~2. Let $f=\sum_na_nq^n$ be the mod~2 modular form associated to~$f$.

As we said already, if $a_2=0$ then Mazur proved that $\m$
satisfies multiplicity one. Gross proves in~\cite{gross} that if $a_2\not=1$
then $\m$ satisfies multiplicity one, and
in Chapter~12 of~\cite{gross} (see the top of page 494) he states
as an open problem whether multiplicity one holds in the remaining case
$a_2=1$.

Edixhoven~\cite{edixhoven} proves that if $\rho_{\m}$ is
ramified at $2$, then $\m$ satisfies multiplicity~1 (Theorem~9.2,
part~3). Buzzard~\cite{buzzard} and \cite{buzzard_appendix} shows
that if $\rho_{\m}$
is unramified at~2 but if the image
of Frobenius at 2 is not contained in the scalars of $\GL_2(\hecke{q} / \m)$
then $\m$ satisfies multiplicity one; see Proposition
2.4 of~\cite{buzzard}.

The main theorem of this paper is that, for $q \in \{431,503,2089\}$ 
(note that all of these are prime), there is a non-Eisenstein
maximal ideal $\m$
of $\hecke{q}$ lying above~2, such that $(\hecke{q})_{\m}$
is not Gorenstein. As a consequence, these $\m$ do not satisfy mod~2
multiplicity one. In fact, the theorems above served as a very useful
guide to where to search for such $\m$.
 
It is proved in Matsumura~\cite{matsumura}, Theorem~21.3, %\footnote{Lloyd: can you give a more precise reference rather than just referring to this very dense book?}
that if a $\ztwo$-algebra is
finitely generated and a complete intersection, then it is Gorenstein. Hence
the $\hecke{q}$ above are not complete intersections. The work of
Wiles~\cite{wiles} and Taylor-Wiles~\cite{taylor-wiles} on Fermat's Last Theorem proves as a byproduct
that certain Hecke algebras are
complete intersections. Hence this paper gives a
bound on how effective Wiles and Taylor's methods can be in characteristic~2.
 
Each of the maximal ideals $\m\subset\hecke{q}$ that we construct
are non-Eisenstein and have the property that $\hecke{q}/\m=\f{2}$,
so let us first consider a general surjective representation
$\rho:{\rm Gal}(\overline{\qq}/\qq) \rightarrow
\GL_2(\f{2})$. For such a representation,
the trace of $\rho$ at any unramified prime $p$
can be computed if one knows the splitting of $p$ in the extension of
$\qq$ of degree $6$ cut out by
$\rho$.  Let $K$ denote this extension.
If $p$ splits completely then ${\Frob}_p$ has order~1
and if $p$ splits into $3$ primes in~$K$, then ${\Frob}_p$
has order~2. In both cases the trace of $\rho({\Frob}_p)$ is~0.
The other possibility is that $p$ splits into two primes, and then
${\Frob}_p$ has order~3 and the trace of $\rho({\Frob}_p)$ is~$1$.

By the theorems of Mazur, Gross, Edixhoven and Buzzard above, if we wish
to find examples of maximal ideals where multiplicity one fails, we could
adopt the following approach: we firstly search for modular elliptic
curves~$E$ of prime conductor~$q$, such that 2 is unramified in $K=\qq(E[2])$,
the field generated over~$\qq$ by the coordinates of the 2-torsion of~$E$.
We then require that 2 splits completely in $K$ and that~$K$ has degree~6
over~$\qq$. Note that these conditions imply that~$E$ has good ordinary
reduction
at~2. For any such elliptic curves that one may find, the maximal ideals of
$\hecke{q}$ associated to $E[2]$
will be maximal ideals not covered by any of the multiplicity
one results above.

We search for curves like this
by using a computer to search pre-compiled tables, such as
Cremona~\cite{cremona-table}, and compute the associated maximal
ideals of $\hecke{q}$. 
For each such maximal ideal, we then explicitly
construct the completion of $(\hecke{q})_{\m}$ as a subring of a direct
sum of finite extensions of $\qq_2$ and check to see whether it is
Gorenstein. Note that a Noetherian local ring is Gorenstein if and only
if its completion is. Hence this procedure will test the Gorenstein-ness
of localisations of $\hecke{q}$ at non-Eisenstein maximal ideals not
covered by the above theorems.

It is slightly surprising to note that after trying only one or two
examples, one finds maximal ideals where multiplicity one fails.
Our discoveries are summarised in the following theorem.
 
\begin{thm}
Assume that $q \in\{431,503,2089\}$. Then there is a maximal ideal
$\m$ of $\hecke{q}$ such that $(\hecke{q})_{\m}$ is not Gorenstein,
and hence for which mod~2 multiplicity one does not hold.  
\end{thm}
 
Note that as a consequence, the ring $\hecke{q}$ is in these
cases, by definition, not Gorenstein. 

The calculations behind this theorem were made possible with the
\HECKE\ package, running on
the \magma\ computer algebra system \cite{magma}. 
\magma\ includes an environment for
specialised number-theoretic calculations, and also incorporates Cremona's 
database of elliptic curves, and the \HECKE\ package, written by William
Stein, implements efficient algorithms to generate spaces of modular forms
of arbitrary level, weight ($ \ge 2$) and character over global and
finite fields. Without this computing package this paper could not have been
written.
 
An application of this paper's results can be found in Section 6 of
Emerton~\cite{emerton}.
Let $X$ be the free
$\zed$-module of divisors supported on the set of singular points 
of the curve $X_0(q)$ in characteristic $q$. 
Theorem~0.5 of \cite{emerton}
shows that $\hecke{q}$ is Gorenstein if and only if~$X$ 
is an invertible $\hecke{q}$-module. The module~$X$ can be explicitly calculated 
quickly using the Mestre-Oesterle method of graphs, from~\cite{mestre}, as implemented in,
e.g., \magma.

The exact sequence of ${\rm Gal}(\overline{K}/K)$-modules on
page~488 of Gross~\cite{gross}, is split as an exact sequence
of Hecke modules, in the Gorenstein case.
Emerton proves that the
analogue of this short exact sequence in the non-Gorenstein case
is never split.   He uses this to prove results about the $\m$-adic Tate
module of $J_0(q)$.
 
In Ribet-Stein~\cite{ribet-stein} the existence of non-Gorenstein Hecke
algebras is discussed in the context of the level optimisation procedure
associated with Serre's conjecture (see section 3.7.1).
 
I would like to thank Kevin Buzzard for suggesting this problem and for
giving me much encouragement, and William Stein for helpful suggestions and
technical support. I would also like to thank the anonymous referee for all
of his or her suggestions, which have improved the paper and its exposition
markedly. 
\section{$\hecke{431}$}
 
\begin{thm}
There is a maximal ideal~$\m$ of $\hecke{431}$ such that
$(\hecke{431})_{\m}$ is not Gorenstein.
\end{thm}

%This section will prove the Theorem above, by considering the Hecke algebra~$\hecke{431}$.

Set $q=431$.
There are two non-isogenous modular elliptic curves $E_1$ and $E_2$ of
conductor $q$, defined over~$\qq$, with corresponding minimal Weierstrass
equations
\begin{eqnarray*}%{displaymath}
&E_1&\negthickspace\negthickspace\negmedspace:\: y^2 + xy = x^3 - 1 \\
%\end{displaymath}
%and
%\begin{displaymath}
&E_2&\negthickspace\negthickspace\negmedspace:\: y^2 + xy + y = x^3 - x^2 - 9x - 8.
\end{eqnarray*}%{displaymath}
The corresponding modular forms $d_1,d_2 \in S_2(\Gamma_0(431))$ have Fourier
coefficients in $\zed$.  One can easily check
that the fields $\qq(E_i[2])$ are isomorphic, and of degree~6 over~$\qq$,
which shows that the two elliptic curves have isomorphic and irreducible
mod~2 Galois representations. One checks that~2 splits
completely in $K=\qq(E_i[2])$.

Let
$$
\rho: \mbox{\rm Gal}(K/\qq) \longrightarrow \mbox{\rm Aut}(E_1[2])
$$
denote the associated mod~2 Galois representation. We claim that
there are precisely four eigenforms in $S_2(\Gamma_0(431))$ giving
rise to $\rho$. To see this, note firstly that $\rho(\Frob_3)$ has
order~3 and hence has trace~1. We compute the characteristic polynomial of
$T_3$ acting on $S_2(\Gamma_0(431))$ and reduce it modulo~2.
We find that the resulting polynomial is of the form $(X-1)^4g(X)$, where
$g(1)\not=0$. This shows already that there are at
most four eigenforms which could give rise to $\rho$. We now compute
the $q$-expansions of these four candidate eigenforms, as elements
of $\Qtwobar[\![q]\!]$. Two of the eigenforms, say $d_1$ and $d_2$,
corresponding to the elliptic
curves $E_1$ and $E_2$ above, of course have coefficients in $\Z$.
The other two, $d_3$ and $d_4$,
are conjugate and are defined over $\Q_2(\sqrt{10})$.
We now check that the reductions to $\F_2$
of the first~72 $q$-expansion coefficients (the Sturm bound) of all four
of these eigenforms are equal. By~\cite{sturm}, this implies
that the four forms themselves are congruent, and hence all give rise to
isomorphic mod~2 Galois representations. By construction, these
representations must all be isomorphic to $\rho$.

Let $\m$ denote the corresponding maximal ideal of $\T_q$. Our goal now
is to compute the completion $\T$ of $(\T_q)_\m$
explicitly. Let $V$ be the 4-dimensional space over $\Qtwobar$
spanned by the four eigenforms giving rise to $\rho$.
The ring $\T$ is the $\Z_2$-subalgebra of $\End(V)$ generated by the Hecke
operators $T_p$ for $p\not=431$, and $U_{431}$.
By the theory of the Sturm bound, this algebra equals the
$\Z_2$-algebra generated by $T_n$ for $n\leq72$. One readily 
computes this algebra. Let $\alpha$ be the coefficient of $q^3$
in $d_3$. Then $d_3$ and $d_4$ are defined over $\Q_2(\alpha)=\Q_2(\sqrt{10})$.
It turns out that $\T$ is isomorphic to
the subring of $\Z_2\oplus\Z_2\oplus\Z_2[\alpha]$
generated as a $\Z_2$-module by $\{(1,1,1),(0,2,0),(0,0,\alpha+1),(0,0,2)\}$.

This is a local ring, with unique maximal ideal $\m$
generated as a $\Z_2$-module by $\{(2,0,0),(0,2,0),(0,0,\alpha+1),(0,0,2)\}$.

Now we claim that there is a reducible
parameter ideal (recall that a parameter ideal is
an ideal that contains a power of 
the maximal ideal).  This will show that $\hecke{}$ is not Gorenstein 
(see~\cite{matsumura}, Theorem~18.1). 
(Note that  $\hecke{}$ is Cohen-Macaulay, as it has Krull dimension one
and has no nonzero nilpotent
elements (see~\cite{matsumura}, Section~17, page~139)).

The ideal $\p$ generated by $(2,2,\alpha+1)$
is a parameter ideal, since it contains~${\m}^2$. We observe that the ideals
$$\p_1=\left((2,0,0),(0,2,\alpha+1)\right)\quad\text{ and }
\quad \p_2=\left((0,2,0),(2,0,\alpha+1)\right)$$
 have intersection exactly $\p$.
Hence $\p$ is a reducible parameter ideal, so $\T$
is not Gorenstein, and as $\T$ is the completion of the Noetherian
local ring $(\T_{431})_{\m}$
we deduce that this ring is also not Gorenstein. Hence mod~2 multiplicity
one fails for~$\m$.
 
\section{$\hecke{503}$}

\begin{thm}
There is a maximal ideal~$\m$ of $\hecke{503}$ such that
$\left(\hecke{503}\right)_{\m}$ is not Gorenstein.
\end{thm}

The argument is similar to that of the previous section, and we merely
sketch it. In fact, the case $q=503$ is slightly technically simpler than
the case of $q=431$ because all four relevant newforms are defined
over $\Z_2$.

There are three isogeny classes of modular elliptic curves of conductor 503.
Let $F_1$, $F_2$ and $F_3$ denote representatives in each class,
and let $f_1$, $f_2$ and $f_3$ denote the corresponding modular forms.
 One checks that the fields~$\qq(F_i[2])$ generated by the 2-torsion of the
three elliptic curves are all isomorphic. If $K$ denotes this extension,
then one checks that $K$ has degree~6 over~$\Q$, and that 2 is unramified
and splits completely in the intergers of~$K$. Let $\rho:\GQ\to\GL_2(\F_2)$
denote the corresponding Galois representation,
and let $\m$ denote the maximal ideal of $\hecke{503}$ corresponding
to the 2-torsion in any of these curves.
An explicit computation of
$S_2(\Gamma_0(503))$ and the eigenvalues mod~2 of
the Hecke operator $T_{11}$ shows that
there can be at most one other eigenform $f_4$ of level~503 giving
rise to $\rho$, and indeed one can check that such a form $f_4$
exists, defined over $\Q_2$ but not $\Q$ (one need only check
congruence for the first~84 coefficients).
The completion of $(\T_{503})_{\m}$ now can be checked to be isomorphic to the
subring of $(\Z_2)^4$ generated as a $\Z_2$-module by the elements

$$\left\{(1,1,1,1),(0,2,0,0),(0,0,2,2),(0,0,0,4)\right\} \subset
(\ztwo)^4.$$

The unique maximal ideal of this ring is generated as a $\Z_2$-module
by the elements
$$\{(2,0,0,0),(0,2,0,0),(0,0,2,2),(0,0,0,4)\}.$$

We see easily that ${(\hecke{503})}_\m$ is Cohen-Macaulay, as it has Krull
dimension 1 and no nilpotent elements.
The ideal $\p$ generated by $(2, 2, 2, 2)$ is a parameter ideal,
as $\m^2 \subseteq \p$. Finally, the ideals
$$\p_1:=\left((0,2,2,2),(2,0,0,0)\right)\:{\rm and} \:
\p_2:=\left((2,0,2,2),(0,2,0,0)\right)$$
have $\p$ as their intersection, hence $\p$ is reducible, and
therefore the localisation $(\hecke{503})_{\m}$ is not Gorenstein 
and hence $\hecke{503}$ is not Gorenstein. Hence as before, mod~2
multiplicity one fails for~$\m$.
 
\section{Other examples}
 
After these initial examples were discovered, William Stein suggested
another way to check directly that multiplicity
one fails in these and other cases, by an explicit computation in the
Jacobian of the relevant modular curve.

Let~$q$ be a prime and let $E_1$ and $E_2$ denote two non-isogenous
new optimal modular elliptic curves (also called strong Weil
curves) of conductor~$q$, viewed as abelian subvarieties
of the abelian variety $J_0(q)$.  
Suppose that the two modular forms corresponding to $E_1$ and $E_2$
are congruent modulo~$2$. Let $\m$ be the corresponding maximal ideal
of $\hecke{q}$ over~$2$.
If mod~$2$ multiplicity one holds for~$\m$, then 
$J_0(q)[\m]=E_1[2]=E_2[2]$, viewed as subsets of $J_0(q)$. In particular, 
$E_i(\C) \bigcap E_j(\C)$ is non-zero.
But this intersection can be explicitly computed using the
{\tt IntersectionGroup} command in \magma, and if it is zero then
we have verified that mod~2 multiplicity one has failed
without having to compute the completions of the relevant Hecke algebra
explicitly.

In the case $q=431$, we find that $E_1(\C) \bigcap E_2(\C) = \{0\}$,
so  we have another proof that mod~2 multiplicity one and the
Gorenstein property fail for the Hecke algebra of level~431.
A similar argument verifies failure of mod~$2$ multiplicity one at
level~503. 

This method can also be used to check our third example.
There are five isogeny classes of elliptic curves with conductor
2089, of which four, the ones labelled 2089A, 2089C, 2089D, 2089E in
Cremona's tables, have ismorphic mod~$2$ representations.  (Note that
2089B, has rational 2-torsion, so its associated Galois representation
is reducible.) Using the results of the previous section, we find that the
intersection in $J_0(2089)$ of the elliptic curves labeled
2089A and 2089B is $\{0\}$.  (Incidentally, the intersection of the curves
labeled 2089A and 2089E is $(\zed/2\zed)^2$.)
Hence we have another example of failure of mod~2 multiplicity one
and Gorenstein-ness.
 
This paper answers the question raised in \cite{buzzard} and in
\cite{gross}, in that it exhibits specific examples of non-Gorenstein Hecke
algebras. This raises the natural question of deciding exactly
which Hecke algebras are not Gorenstein, without explicitly computing them,
and proving theorems like those of Buzzard, Edixhoven, Gross and Mazur for
these algebras. A natural next step is to ask the following:

\begin{question}
Are there infinitely many prime integers~$q$ such that
$\hecke{q}$ is not Gorenstein, and hence where mod~2 multiplicity one fails?
\end{question}

L.\ J.\ P.\ Kilford,

Imperial College,

Department of Mathematics,

180 Queen's Gate,

London,

SW7 2AZ,

United Kingdom.

\medskip

{\tt l.kilford@ic.ac.uk}

\medskip

(landline) 01502 566564.
%(mobile) 07752 340186.
(fax) 0207 5948517.

\end{document}